\documentclass[journal]{IEEEtran}

\ifCLASSINFOpdf

\else

\fi
\usepackage{graphicx}          

\usepackage[dvips]{epsfig}    
\usepackage{latexsym,amssymb}
\usepackage{amsmath, amsbsy, amsthm}
\usepackage{amsopn, amstext}
\usepackage{colortbl}
\usepackage{cite}
\newtheorem{remark}{Remark}
\newtheorem{theorem}{Theorem}

\newtheorem{lemma}{Lemma}

\newtheorem{problem}{Problem}
\newtheorem{assumption}{Assumption}
\newtheorem{corollary}{Corollary}
\allowdisplaybreaks
\begin{document}

\title{The Difference and Unity of Irregular LQ Control and Standard LQ Control and Its Solution}

\author{Huanshui~Zhang~and~Juanjuan~Xu~
\thanks{*This work is supported by the National Natural Science Foundation of China under Grants 61633014, 61873332, U1806204, U1701264, 61922051, the foundation for Innovative Research Groups of National Natural Science Foundation of China (61821004) and Youth Innovation Group Project of Shandong University (2020QNQT016).}
\thanks{H. Zhang is with Shandong University and Shandong University of
Science and Technology, Shandong, P.R.China (most of the work was completed in Shandong University). 
Juanjuan Xu is with School of Control Science and Engineering of Shandong University, Shandong, P.R.China.
        {\tt\small hszhang@sdu.edu.cn, juanjuanxu@sdu.edu.cn}}

}

\maketitle

\begin{abstract}
Irregular linear quadratic control (LQ, was called Singular LQ) has been a long-standing problem since 1970s.
This paper will show that an irregular LQ control (deterministic) is solvable (for arbitrary initial value) if and only if
the LQ cost can be rewritten as a regular one by changing the terminal cost $x'(T)Hx(T)$ to $x'(T)[H+P_1(T)]x(T)$, while the optimal controller can achieve $P_1(T)x(T)=0$ at the same time.
In other words, the irregular controller (if exists) needs to do two things at the same time, one thing is to minimize the cost and the other is to achieve the terminal constraint $P_1(T)x(T)=0$, which clarifies the essential difference of irregular LQ from the standard LQ control where the controller is to minimize the cost only.

With this breakthrough, we further study the irregular LQ control for stochastic systems with multiplicative noise. A sufficient solving condition and the optimal controller is presented based on Riccati equations.

\end{abstract}


%
\IEEEpeerreviewmaketitle

\section{Introduction}

~~~~Linear-quadratic (LQ) optimal control has received much attention in recent years due to the widely applications in modern engineering \cite{Bell,Bellman,Kalman}.
Considering the singularity of the weighting matrix of the control in the cost function, LQ optimal control problem is mainly consisting of regular optimal control
and irregular optimal control. Most of the previous works have been focused on regular case. In particular, when the weighting matrix of the control in the cost function is positive-definite, the LQ optimal control naturally belongs to
the regular case which has been extensively studied in \cite{Kalman,Letov, mp,bellman}.
When the weighting matrix of the control in the cost function is in more general case of indefinition, \cite{zhouxunyu} studied the stochastic optimal control and obtained the optimal solution where the stochastic Riccati equation is strictly required to be regular, i.e., the results are only applicable to regular LQ problems.

In the case of irregularity, the optimal LQ control has been remaining major challenging although much efforts have been made since 1970's. In \cite{Gurman,Moore,Williems} and references therein, the singular LQ control was studied by using `Transformation in state space', where the problem with control weighting matrix is zero ($R=0$) was studied. It was shown that the problem is solvable if the initial value is given like $x_2(0)=C_{21}(0)x_1(0)$. Otherwise, an impulse control must be applied at the initial time \cite{Ho}. In other words, the approach of `Transformation in state space' is only applicable to the case of specified initial value. In \cite{Krener,Bonnans,XuZhang}, the approach of `higher order maximum principle' was applied to singular LQ control. However, if the higher derivatives vanish, it is impossible to find the singular control with this approach \cite{Gabasov}. The third approach is the perturbation approach in \cite{chenhanfu}, \cite{sunliyong}. The optimal solution is obtained by using the limitation of the solution to Riccati equation when the perturbation is approaching to zero.

More recently, with the analytical solution to a forward and backward differential equations (FBDEs), \cite{zhangIR} considered the irregular LQ control for deterministic systems with arbitrary initial value, where the irregular controller was designed based on a regular Riccati equation and the controllability of a subsystem.

In this paper we will study the irregular optimal control problem aiming to explore the difference between the regular control and irregular control (see Theorem \ref{thm} in the below). It is interesting to show that an irregular LQ control is solvable if and only if LQ cost can be rewritten as a regular one by changing the terminal term of the LQ cost, while the controller can make the changed terminal to be zero. Moreover, we extend the results to stochastic control problem with irregular cost (Theorem \ref{thsto}).

 %

The remainder of the paper is organized as follows. Section II presents the solution for the deterministic optimal control problem with irregular performance.
The solution to the stochastic optimal control problem with irregular performance is given in Section III.
Some concluding words are given in Section IV. Some proofs of the results are presented in Appendix.

The following notations will be used throughout this paper: $R^n$
denotes the family of $n$ dimensional vectors. $x'$ means the
transpose of $x.$ It is defined that $\|x\|^2=x'x.$ A symmetric matrix $M>0\ (\geq 0)$ means
strictly positive definite (positive semi-definite). $Range(M)$ represents the range of the matrix $M.$
$M^{\dag}$ is called  the Moore-Penrose inverse \cite{pinv} of the matrix $M$ if it satisfies $MM^{\dag}M=M,~M^{\dag}MM^{\dag}=M^{\dag},
(MM^{\dag})¡¯=MM^{\dag}$ and $(M^{\dag}M)'=M^{\dag}M.$


\section{Deterministic optimal control with irregular performance}

In this section, we consider the deterministic optimal control with irregular performance where the
linear system governed by a differential equation:
\begin{eqnarray}
\dot{x}(t)&=&A(t)x(t)+B(t)u(t),~x(t_0)=x_0,\label{d1}
\end{eqnarray}
where $x\in R^n$ is the state, $u\in R^m$ is the control input.
The matrices $A,B,\bar{A},\bar{B}$ are constant matrices with appropriate dimension. $x_0$ represents the initial value.
The cost function is given by
\begin{eqnarray}
J_0(t_0,x_0;u)&=&\int_{t_0}^T[x'(t)Q(t)x(t)+u'(t)R(t)u(t)]dt\nonumber\\
&&+x'(T)Hx(T),\label{d2}
\end{eqnarray}
where $Q(t)\geq 0,R(t)\geq 0$ are symmetric matrices with appropriate dimensions.
\begin{problem}
For any $(t_0,x_0)$, find a controller $u(t)$ such that (\ref{d2}) is minimized subject to (\ref{d1}).
\end{problem}
Noting that $R(t)$ is semi positive-definite, the problem was usually called singular control \cite{Gurman,Williems}, which remains to be solved due to the difficulty caused by the {\em regularity}.

\subsection{What is an irregular LQ control}
Singular LQ control contains regular and irregular control,  the first case is easily done using the standard control approach, the second case of irregular is much involved as said in the above. To define irregular control problem, we introduce the following Riccati equation associated with system (\ref{d1}) and cost (\ref{d2})
\begin{eqnarray}
0&=&\dot{P}(t)+A'(t)P(t)+P(t)A(t)+Q(t)\nonumber\\
&&-P(t)B(t)R^{\dag}(t)B'(t)P(t), \label{5}
\end{eqnarray}
where the terminal condition is given by $P(T)=H$ and $R^{\dag}(t)$ represents the Moore-Penrose inverse of $R(t)$.  If $Range [B'(t)P(t)]\subseteq Range [R(t)]$, the LQ control problem is standard and called regular.  Otherwise if
\begin{eqnarray}
Range [B'(t)P(t)]\not\subseteq Range [R(t)],\label{d33}
\end{eqnarray}
the LQ control is called irregular and the performance cost (\ref{d2}) is irregular accordingly.

\subsection{Why is it difficult?}
The irregularity implies the controller is unsolvable with classical control theory.  In fact,
the irregularity leads to extremely difficulty to obtain the controller. To show this, we present an example where the system is governed by
\begin{eqnarray}
\dot{x}(t)&=&x(t)+\left[
                            \begin{array}{cc}
                              1 & -1 \\
                            \end{array}
                          \right]u(t),~x(t_0)=x_0,\label{de2}
\end{eqnarray} and the cost function is given by \begin{eqnarray}
J_T(t_0,x_0;u)=\int_{t_0}^Tu'(t)\left[
                                   \begin{array}{cc}
                                     1 & 0 \\
                                     0 & 0 \\
                                   \end{array}
                                 \right]u(t)dt+x'(T)x(T).\label{de1}
\end{eqnarray}
The solution to the Riccati equation $0=\dot{P}(t)+2P(t)-P^2(t)$ with $P(T)=1$ is given by $P(t)=\frac{2}{1+e^{2(t-T)}}.$
Then, it holds that $Range [B'P(t)]\nsubseteq Range (R)$ where $B=\left[
                            \begin{array}{cc}
                              1 & -1 \\
                            \end{array}
                          \right]$ and $R=\left[
                                   \begin{array}{cc}
                                     1 & 0 \\
                                     0 & 0 \\
                                   \end{array}
                                 \right]$. This implies that it is unable to obtain $u(t)$ from the classical equilibrium condition $Ru(t)+B'P(t)x(t)=0$ for arbitrary $x(t)$.\\

The irregularity also leads to fundamental difficulty, that is, completing sum of squares can not be achieved for irregular LQ cost (\ref{d2}). Actually,
for the above optimization problem of minimizing (\ref{de1}) subject to (\ref{de2}), by taking derivative to $x'(t)P(t)x(t),$ it yields that
\begin{eqnarray}
\frac{d}{dt}[x'(t)P(t)x(t)]=2u'(t)B'P(t)x(t)+x'(t)P^2(t)x(t).\nonumber
\end{eqnarray}
Then the cost function (\ref{de1}) can be rewritten by taking integration from $t_0$ to $T$ in the above equation as
\begin{eqnarray}
&&J_T(t_0,x_0;u)\nonumber\\
&=&x'(t_0)P(t_0)x(t_0)+\int_{t_0}^T\Big[u'(t)Ru(t)+2u'(t)\nonumber\\
&&\times B'P(t)x(t)+x'(t)P^2(t)x(t)\Big]dt\nonumber\\
&=&x'(t_0)P(t_0)x(t_0)+\int_{t_0}^T\Big\{\Big[u(t)+R^{\dag}B'P(t)x(t)\Big]'\nonumber\\
&&\times R\Big[u(t)+R^{\dag}B'P(t)x(t)\Big]Ru(t)\nonumber\\
&&+2u'(t)(I-RR^{\dag})B'P(t)x(t)\Big\}dt.\nonumber
\end{eqnarray}
From the above cost function and the fact that $Range [B'P(t)]\nsubseteq Range (R),$ it is seen that the optimal controller can not be obtained by completing sum of squares
because the last term in the above is not zero.

Thus the irregularity leads to the invalidity of standard methods for LQ control. In this paper, in order to find a way to solve the irregular control,
we will first explore the essential difference of irregular control from regular one in the following Theorem 1.

\subsection{Solution to deterministic optimal control with irregular performance}

Firstly, we present the maximum principle for Problem 1 \cite{zhangIR}.

\begin{lemma}\label{lem4}
If Problem 1 is solvable, then the optimal controller satisfies
\begin{eqnarray}
0&=&R(t)u(t)+B'(t)p(t),\label{8}
\end{eqnarray}
where $p(t)$ obeys the following dynamics:
\begin{eqnarray}
\dot{p}(t)&=&-A'(t)p(t)-Q(t)x(t),\label{4}
\end{eqnarray}
with the terminal value $p(T)=Hx(T)$. Conversely, if FBDEs (\ref{d1}), (\ref{4}) and (\ref{8}) is solvable, then Problem 1 is also solvable.
\end{lemma}
\emph{Proof.} The proof follows from the maximum principle by using the fact that $Q(t)\geq0, R(t)\geq 0$. So we omit the details. \hfill $\blacksquare$

We next make some denotations for convenience of the derivation of the main result.
Let $rank (R(t))=m_0< m,$ thus $rank\Big[I-R^{\dag}(t)R(t)\Big]=m-m_0>0$. There is an elementary row transformation matrix $T_0(t)$ such that
\begin{eqnarray}
T_0(t)\Big[I-R^{\dag}(t)R(t)\Big]=\left[
                              \begin{array}{c}
                                 0  \\
                                 \Upsilon_{T_0}(t)\\
                              \end{array}
                            \right], \label{jnYC1}
\end{eqnarray}
where $\Upsilon_{T_0}(t)\in R^{[m-m_0]\times m}$ is full row rank. Furthermore denote
\begin{eqnarray}
A_0(t)&=&A(t)-B(t)R^{\dag}(t)B'(t)P(t),\nonumber\\
D_0(t)&=&-B(t)R^{\dag}(t)B'(t),\nonumber\\
\left[ \begin{array}{cc}\ast & B_0(t)\\
                              \end{array}
                            \right] &=&B(t)\Big[I-R^{\dag}(t)R(t)\Big]{T_0}^{-1}(t),\nonumber\\
\left[
  \begin{array}{cc}
    * & G_0(t) \\
  \end{array}
\right]&=&T_0^{-1}(t),\nonumber
\end{eqnarray}
and define
\begin{eqnarray}
0&=&\dot{P}_1(t)+P_1(t)A_0(t)+A_0'(t)P_1(t)\nonumber\\
&&+P_1(t)D_0(t)P_1(t), \label{d3}
\end{eqnarray}
where the terminal value $ P_1(T)$ is to be determined.

We are now in the position to give the main result of this section as follows.

\begin{theorem}\label{thm}
Problem 1 is solvable if and only if there exists a matrix $P_1(T)$ satisfying $0=B_0'(T)[P(T)+P_1(T)]$
such that the following changed cost
\begin{eqnarray}
{\bar J}_0(t_0,x_0;u)&=&J_0(t_0,x_0;u)+x'(T)P_1(T)x(T)\label{x1}
\end{eqnarray}
is regular and $P_1(T)x(T)=0$ is achieved with the controller minimizing (\ref{x1}).

\end{theorem}

{\em Proof}. ``Sufficiency" The aim is to prove if there exists a matrix $P_1(T)$ satisfying $0=B_0'(T)[P(T)+P_1(T)]$
such that (\ref{x1}) is regular and $P_1(T)x(T)=0$ is achieved, then Problem 1 is solvable. Based on Lemma \ref{lem4},
it is sufficient to show that the FBDEs (\ref{d1}), (\ref{4}) and (\ref{8}) is solvable.  To this end, we will verify that the following $(p(t),x(t))$ solves the FBDEs:
\begin{eqnarray}
p(t)&=&P(t)x(t)+P_1(t)x(t),\label{13}\\
\dot{x}(t)&=&\Big\{A(t)-B(t)R^{\dag}(t)B'(t)\Big[P(t)+P_1(t)\Big]\Big\}x(t)\nonumber\\
&&+B(t)\Big[I-R^{\dag}(t)R(t)\Big]z(t),\label{1}
\end{eqnarray}
where $P(t)$ is defined in (\ref{5}), $P_1(t)$ is defined by (\ref{d3}) with $P_1(T)$ satisfying $0=B_0'(T)[P(T)+P_1(T)]$ and $z(t)$ is a vector with
compatible dimension such that $P_1(T)x(T)=0$.

The verification is divided into two steps. The first step is to prove that the Riccati equation $P(t)+P_1(t)$ is regular.
From the regularity of (\ref{x1}), we have the following Riccati equation is regular:
\begin{eqnarray}
0&=&\dot{\bar{P}}(t)+A'(t)\bar{P}(t)+\bar{P}(t)A(t)+Q(t)\nonumber\\
&&-\bar{P}(t)B(t)R^{\dag}(t)B'(t)\bar{P}(t),\label{dr1}
\end{eqnarray}
with terminal value $\bar{P}(T)=H+P_1(T)$. That is,
\begin{eqnarray}
\Big[I-R(t)R^{\dag}(t)\Big]B'(t)\bar{P}(t)=0.\label{12}
\end{eqnarray}
In addition, note that $P(t)+P_1(t)$ satisfies
\begin{eqnarray}
0&=&\dot{P}(t)+\dot{P}_1(t)+A'(t)\Big[P(t)+P_1(t)\Big]\nonumber\\
&&+\Big[P(t)+P_1(t)\Big]A(t)+Q(t)-\Big[P(t)+P_1(t)\Big]\nonumber\\
&&\times B(t)R^{\dag}(t)B'(t)\Big[P(t)+P_1(t)\Big], \label{7}
\end{eqnarray}
with the same terminal value $\bar{P}(T)=H+P_1(T)$ to (\ref{dr1}).
This implies that
\begin{eqnarray}
\bar{P}(t)=P(t)+P_1(t).\label{d4}
\end{eqnarray}
Thus, it is obtained from (\ref{12}) that
\begin{eqnarray}
&&\Big[I-R(t)R^{\dag}(t)\Big]B'(t)\Big[P(t)+P_1(t)\Big]\nonumber\\
&=&\Big[I-R(t)R^{\dag}(t)\Big]B'(t)\bar{P}(t)=0.\nonumber
\end{eqnarray}

The second step is to take derivatives on the right of (\ref{13}). By using (\ref{d4}), (\ref{dr1}) and (\ref{1}), we derive from (\ref{13}) that
\begin{eqnarray}
\frac{d}{dt}\Big[\bar{P}(t)x(t)\Big]
&=&\frac{d}{dt}\Big[[P(t)+P_1(t)]x(t)\Big]\nonumber\\
&=&-\Big[A'(t)\bar{P}(t)+\bar{P}(t)A(t)+Q(t)\nonumber\\
&&-\bar{P}(t)B(t)R^{\dag}(t)B'(t)\bar{P}(t)\Big]x(t)\nonumber\\
&&+\bar{P}(t)\Big[A(t)-B(t)R^{\dag}(t)B'(t)\bar{P}(t)\Big]x(t)\nonumber\\
&&+\bar{P}(t)B(t)\Big[I-R^{\dag}(t)R(t)\Big]z(t)\nonumber\\
&=&-A'(t)\bar{P}(t)x(t)-Q(t)x(t),\nonumber
\end{eqnarray}
where the last term in the first equality is equal to zero by using (\ref{12}).
Furthermore, by denoting $u(t)=-R^{\dag}(t)\bar{P}(t)x(t)+[I-R(t)R^{\dag}(t)]z(t)$, it follows that $0=R(t)u(t)+B'(t)\bar{P}(t)x(t)$.
In additionally, (\ref{1}) can be rewritten as $\dot{x}(t)=A(t)x(t)+B(t)u(t)$. Accordingly, $(p(t), x(t))$ defined by (\ref{13}) and (\ref{1})
solves the FBDEs (\ref{d1}), (\ref{4}) and (\ref{8}). Thus, based on Lemma \ref{lem4}, Problem 1 is solvable.

``Necessity" The aim is to verify the regularity of (\ref{x1}) and $P_1(T)x(T)=0$. In fact, by using Theorem 2 in \cite{zhangIR}, the necessary
condition for the solvability of Problem 1 is that $B_0'(t)[P(t)+P_1(t)]=0$ and $P_1(T)x(T)=0$ holds.
Thus, the key is to prove (\ref{x1}) is regular, that is, (\ref{12}) holds where $\bar{P}(t)$ satisfies (\ref{dr1}).
In fact, by using (\ref{d4}) and the fact that $\Upsilon_{T_0}(t)$ has full row rank, it is obtained from $B_0'(t)[P(t)+P_1(t)]=0$ that
\begin{eqnarray}
0&=&\Upsilon_{T_0}'(t)B_0'(t)\bar{P}(t)\nonumber\\
&=&\left[
     \begin{array}{cc}
       0 & \Upsilon_{T_0}'(t) \\
     \end{array}
   \right]\left[
            \begin{array}{c}
              * \\
              B_0'(t) \\
            \end{array}
          \right]\bar{P}(t)\nonumber\\
&=&\Big[I-R(t)R^{\dag}(t)\Big]T_0'(t)[T_0'(t)]^{-1}\nonumber\\
&&\times \Big[I-R(t)R^{\dag}(t)\Big]B'(t)\bar{P}(t)\nonumber\\
&=&\Big[I-R(t)R^{\dag}(t)\Big]B'(t)\bar{P}(t).\nonumber
\end{eqnarray}
This implies that (\ref{x1}) is regular. The proof is now completed. \hfill $\blacksquare$

\begin{remark}
\begin{itemize}

\item It is obvious that $P_1(T)=0$ for the regular (standard) LQ control, while $P_1(T)\not=0$ for the irregular LQ control. So an essential difference of irregular LQ from regular one is explored that the irregular controller (if exists) needs to do two things at the same time, one thing is to minimize the cost (\ref{x1}) and the other is to achieve the terminal constraint $P_1(T)x(T)=0$.
\item Though the difference, the LQ control problem (irregular and regular) can be solved in a unified way as in Theorem \ref{thm}.\end{itemize}

\end{remark}

To conclude this section, we present the optimal controller of Problem 1.
\begin{corollary}\label{cor1}
If there exists a matrix $P_1(T)$ satisfying $0=B_0'(T)[P(T)+P_1(T)]$
such that the changed cost (\ref{x1}) is regular, then the optimal controller is given by
\begin{eqnarray}
u(t)&=&-R^{\dag}(t)B'(t)\Big[P(t)+P_1(t)\Big]x(t)\nonumber\\
&&+G_0(t)u_1(t),\label{x3}
\end{eqnarray}
where $P(t)$ and $P_1(t)$ satisfy Riccati equations (\ref{5}) and (\ref{d3}), and $u_1(t)$ is chosen such that $P_1(T)x(T)=0$.
The optimal cost is given by
\begin{eqnarray}
J^*(t_0,x_0;u)=x_0'\Big[P(t_0)+P_1(t_0)\Big]x_0.\label{6}
\end{eqnarray}
\end{corollary}
\emph{Proof.} By solving the regular optimal control problem of minimizing (\ref{x1}) subject to (\ref{d1}), we have the optimal control
is given by
\begin{eqnarray}
u(t)&=&-R^{\dag}(t)B'(t)\bar{P}(t)x(t)\nonumber\\
&&+\Big[I-R^{\dag}(t)R(t)\Big]z(t),\label{d5}
\end{eqnarray}
where $z(t)$ is chosen such that $P_1(T)x(T)=0$. Combining with the denotation above (\ref{d3}), we can rewrite the last term in the above equation as
\begin{eqnarray}
\Big[I-R^{\dag}(t)R(t)\Big]z(t)&=&T_0^{-1}(t)T_0(t)\Big[I-R^{\dag}(t)R(t)\Big]z(t)\nonumber\\
&=&\left[
  \begin{array}{cc}
    * & G_0(t) \\
  \end{array}
\right]\left[
                                \begin{array}{c}
                                  0 \\
                                  \Upsilon_{T_0}(t) \\
                                \end{array}
                              \right]z(t)\nonumber\\
&=&G_0(t)\Upsilon_{T_0}(t)z(t).\nonumber
\end{eqnarray}
By letting $\Upsilon_{T_0}(t)z(t)=u_1(t)$, the optimal controller (\ref{x3}) follows. The proof is now completed. \hfill $\blacksquare$

\section{Stochastic optimal control with irregular performance}

In this section,  we will extend the above results to the stochastic optimal control with irregular performance where the linear control system
is governed by  It\^{o} stochastic differential equation:
\begin{eqnarray}
dx(t)&=&\Big[A(t)x(t)+B(t)u(t)\Big]dt+\Big[\bar{A}(t)x(t)\nonumber\\
&&+\bar{B}(t)u(t)\Big]dw(t),~x(t_0)=x_0,\label{s1}
\end{eqnarray}
where $x\in R^n$ is the state, $u\in R^m$ is the control input. $w(t)$ is a standard one-dimension Brownian motion.
The filtration $\mathcal{F}_t$ is generated by $\{w(t), t\geq t_0\}$, that is, $\mathcal{F}_t=\sigma\{w(s), t_0\leq s\leq t\}$.
The matrices $A(t),B(t),\bar{A}(t),\bar{B}(t)$ are deterministic matrices with compatible dimensions.
The cost function is given by
\begin{eqnarray}
J(t_0,x_0;u)&=&E\int_{t_0}^T\Big[x'(t)Q(t)x(t)+u'(t)R(t)u(t)\Big]dt\nonumber\\
&&+Ex'(T)Hx(T),\label{s2}
\end{eqnarray}
where $Q(t),R(t),H$ are symmetric matrices with compatible dimensions. The set of the admissible controllers is denoted by
\begin{eqnarray}
\mathcal{U}[t_0,T]&=&\Big\{u(t),t\in[t_0,T]\Big|u(t)~\mbox{is}~\mathcal{F}_{t}~\mbox{adapted},\nonumber\\
&&
E\int_{t_0}^T\|u(t)\|^2dt<\infty\Big\}.\nonumber
\end{eqnarray}

\begin{problem}
For any $(t_0,x_0)$, find an ${\cal F}_t$-adapted controller $u(t)$ such that (\ref{s2}) is minimized subject to (\ref{s1}).
\end{problem}

To guarantee the solvability of Problem 2, we make the following assumption:
\begin{assumption}\label{assum-con} Convexity
\begin{eqnarray}
J(t_0,0;u)\geq0.\nonumber
\end{eqnarray}
\end{assumption}

Under Assumption \ref{assum-con}, we have the following maximum principle for Problem 2 \cite{hszhang1}.

\begin{lemma}\label{lem1}
If Problem 2 is solvable, then the optimal controller satisfies
\begin{eqnarray}
0=R(t)u(t)+B'(t)p(t)+\bar{B}'(t)q(t),\label{3}
\end{eqnarray}
where $(p(t),q(t))$ obey a backward stochastic differential equation (BSDE):
\begin{eqnarray}
dp(t)&=&-[A'(t)p(t)+\bar{A}'q(t)+Q(t)x(t)]dt\nonumber\\
&&+q(t)dw(t),\label{2}
\end{eqnarray}
with the terminal value as $p(T)=Hx(T)$. Conversely, if FBSDEs (\ref{s1}), (\ref{2}) and (\ref{3}) is solvable, then Problem 2 is also solvable.
\end{lemma}
\emph{Proof.} The proof can be found in Lemma 1 in \cite{hszhang1}. \hfill $\blacksquare$

Parallel to (\ref{5}), we introduce the generalized Riccati equation:
\begin{eqnarray}
0&=&\dot{P}(t)+A'(t)P(t)+\bar{A}'(t)P(t)\bar{A}(t)+P(t)A(t)\nonumber\\
&&+Q(t)-\Gamma_0'(t)\Upsilon_0^{\dag}(t)\Gamma_0(t), \label{r1}
\end{eqnarray}
where
\begin{eqnarray}
\Upsilon_0(t)&=&R(t)+\bar{B}'(t)P(t)\bar{B}(t),\label{Intr3}\\
\Gamma_0(t)&=&B'(t)P(t)+\bar{B}'(t)P(t)\bar{A}(t),\label{Intr4}
\end{eqnarray}
and the terminal condition is given by $P(T)=H$.

As has been studied in Section II, we will focus on the stochastic optimal control with irregular performance, that is,
\begin{eqnarray}
Range [\Gamma_0(t)]\not\subseteq Range [\Upsilon_0(t)].\label{s3}
\end{eqnarray}

\subsection{Preliminaries on stochastic optimal control with irregular performance}

In this subsection, we firstly make some denotations for convenience of use. Without loss of generality, we assume
that $rank \Big[\Upsilon_0(t)\Big]=m_0< m$. Thus $rank\Big[I-\Upsilon_0^{\dag}(t)\Upsilon_0(t)\Big]=m-m_0 >0$. It is not difficult
to know that there is an elementary row transformation matrix $T_0(t)$ such that
\begin{eqnarray}
T_0(t)\Big[I-\Upsilon_0^{\dag}(t)\Upsilon_0(t)\Big]=\left[
                              \begin{array}{c}
                                 0  \\
                                 \Upsilon_{T_0}(t)\\
                              \end{array}
                            \right], \label{jnYC1}
\end{eqnarray}
where $\Upsilon_{T_0}(t)\in R^{(m-m_0)\times m}$ is full row rank. Furthermore, we make the following denotations:
\begin{eqnarray}
 \left[ \begin{array}{cc} \ast & C_0'(t) \\
                              \end{array} \right] &=&\Gamma_0'(t)\Big[I-\Upsilon^{\dag}_0(t)\Upsilon_0(t)\Big]{T_0}^{-1}(t),\nonumber\\
\left[ \begin{array}{cc}\ast & B_0(t)\\
                              \end{array}
                            \right] &=&B(t)\Big[I-\Upsilon_0^{\dag}(t)\Upsilon_0(t)\Big]{T_0}^{-1}(t),\nonumber\\
\left[ \begin{array}{cc}\ast & \bar{B}_0(t)\\
                              \end{array}
                            \right] &=& \bar{B}(t)\Big[I-\Upsilon_0^{\dag}(t)\Upsilon_0(t)\Big] {T_0}^{-1}(t),\nonumber\\
\left[
  \begin{array}{cc}
    * & G(t) \\
  \end{array}
\right]&=&T_0^{-1}(t),\nonumber\\
A_0(t)&=&A(t)-B(t)\Upsilon_0^{\dag}(t)\Gamma_{0}(t),\nonumber\\
\bar{A}_0(t)&=&\bar{A}(t)-\bar{B}(t)\Upsilon_0^{\dag}(t)\Gamma_{0}(t),\nonumber\\
D_0(t)&=&-B(t)\Upsilon_0^{\dag}(t)B'(t),\nonumber\\
\bar{D}_0(t)&=&-\bar{B}(t)\Upsilon_0^{\dag}(t)B'(t),\nonumber\\
F_0(t)&=&-B(t)\Upsilon_0^{\dag}(t)\bar{B}'(t),\nonumber\\
\bar{F}_0(t)&=&-\bar{B}(t)\Upsilon_0^{\dag}(t)\bar{B}'(t),\nonumber
\end{eqnarray}
where $C_0'(t),B_0(t),\bar{B}_0(t)\in R^{n\times (m-m_0)}$, $G(t)\in R^{m\times (m-m_0)}$ and define the following Riccati equation:
\begin{eqnarray}
0&=&\dot{P}_1(t)+P_1(t)A_0(t)+A_0'(t)P_1(t)+P_1(t)D_0(t)P_1(t)\nonumber\\
&&+\Big[\bar{A}'_0(t)+P_1(t)F_0(t)\Big]\Big[I-P_1(t)\bar{F}_{0}(t)\Big]^{\dag}\nonumber\\
&&\times P_1(t)\Big[\bar{A}_0(t)
+\bar{D}_0(t)P_1(t)\Big],\label{p1}
\end{eqnarray}
where the terminal value $P_1(T)$ is to be determined.
Moreover, we define that
\begin{eqnarray}
A_1(t)&=&A_0(t)+D_0(t)P_1(t)+F_0(t)\Big[I-P_1(t)\bar{F}_{0}(t)\Big]^{\dag}\nonumber\\
&&\times P_1(t)\Big[\bar{A}_0(t)+\bar{D}_0(t)P_1(t)\Big],\nonumber\\
B_1(t)&=&B_{0}(t)+F_0(t)\Big[I-P_1(t)\bar{F}_{0}(t)\Big]^{\dag}P_1(t)\bar{B}_0(t),\nonumber\\
\bar{A}_1(t)&=&\bar{A}_0(t)+\bar{D}_0(t)P_1(t)+\bar{F}_0(t)\Big[I-P_1(t)\bar{F}_{0}(t)\Big]^{\dag}\nonumber\\
&&\times P_1(t)\Big[\bar{A}_0(t)+\bar{D}_0(t)P_1(t)\Big],\nonumber\\
\bar{B}_1(t)&=&\bar{B}_{0}(t)+\bar{F}_0(t)\Big[I-P_1(t)\bar{F}_{0}(t)\Big]^{\dag}P_1(t)\bar{B}_0(t).\nonumber
\end{eqnarray}

In view of the Riccati equations (\ref{s1}) and (\ref{p1}), we make the following assumptions for the solutions $P(t)$ and $P_1(t)$ of (\ref{s1}) and (\ref{p1}).

\begin{assumption}\label{assum2}

\begin{enumerate}
  \item \begin{eqnarray}
L'(t)&=&L'(t)\Big[I-P_1(t)\bar{F}_{0}(t)\Big]^{\dag}\Big[I-P_1(t)\bar{F}_{0}(t)\Big],\nonumber \\ \label{p11}
\end{eqnarray}
where $L(t)$ may be $B(t)$, $\bar{A}(t)$, or $\bar{B}(t)$.
  \item \begin{eqnarray}
&&\Big[\Upsilon_0(t)+\bar{B}'(t)P_1(t)\bar{B}(t)\Big]^{\dag}L'(t)\nonumber\\
&=&\Big\{I-\Upsilon_0^{\dag}(t)\bar{B}'(t)\Big[I-P_1(t)\bar{F}_{0}(t)\Big]^{\dag}\nonumber\\
&&\times P_1(t)\bar{B}(t)\Big\}\Upsilon_0^{\dag}(t)L'(t),\label{p2}
\end{eqnarray}
where $L(t)$ may be $B(t)$, $\bar{A}(t)$, or $\bar{B}(t)$.
  \item \begin{eqnarray}
0&=&\bar{B}_{0}'(t)\Big[I-P_1(t)\bar{F}_{0}(t)\Big]^{\dag}P_1(t)\bar{B}_{0}(t). \label{pz1}
\end{eqnarray}

\end{enumerate}
\end{assumption}

\begin{remark}
It is noted that the above assumption is not restrictive. In fact, conditions (\ref{p11}) and (\ref{p2}) hold when the Moore-Penrose inverse become inverse. Also, (\ref{p11}), (\ref{p2}) and (\ref{pz1})
hold naturally for deterministic systems.
\end{remark}

Based on the assumption, we present the following lemmas which is useful for the derivation of the main result.

\begin{lemma}\label{leminv}
Under the assumption (\ref{p11}), it holds that
\begin{enumerate}
  \item  Commutative law
  \begin{eqnarray}
  &&L_1'(t)P_1(t)\Big[I-\bar{F}_0(t)P_1(t)\Big]^{\dag}L_2(t)\nonumber\\
  &=&L_{1}'(t)\Big[I-P_1(t)\bar{F}_0(t)\Big]^{\dag}P_1(t)L_{2}(t),\nonumber
  \end{eqnarray}
  where $L_1(t),L_2(t)$ may be $B(t), \bar{A}(t)$, or $\bar{B}(t)$.
  \item Formula of Moore-Penrose inverse for sum of matrices
  \begin{eqnarray}
  &&\Big[I-\bar{F}_0(t)P_1(t)\Big]^{\dag}L(t)\nonumber\\
  &=&\Big\{I+\bar{F}_0(t)\Big[I-P_1(t)\bar{F}_{0}(t)\Big]^{\dag}P_1(t)\Big\}L(t),\nonumber
  \end{eqnarray} where $L(t)$ may be $B(t), \bar{A}(t)$, or $\bar{B}(t)$.

\end{enumerate}
\end{lemma}
\emph{Proof.} The proof is given in Appendix \ref{a2}. \hfill $\blacksquare$

By using Lemma \ref{leminv}, we have a uniform equation for $P(t)+P_1(t)$ as follows.

\begin{lemma}\label{lem3}
Under the assumption (\ref{p11})-(\ref{p2}), it holds that $\bar{P}(t)=P(t)+P_1(t)$ satisfies the following Riccati equation:
\begin{eqnarray}
0&=&\dot{\bar{P}}(t)+A'(t)\bar{P}(t)+\bar{P}(t)A(t)+\bar{A}'(t)\bar{P}(t)\bar{A}(t)\nonumber\\
&&+Q(t)-\bar{\Gamma}'(t)\bar{\Upsilon}^{\dag}(t)\bar{\Gamma}(t),\label{s7}
\end{eqnarray}
where $P(t)$ and $P_1(t)$ are solutions of (\ref{r1}) and (\ref{p1}) respectively, the terminal value is given by $\bar{P}(T)=H+P_1(T)$ and
\begin{eqnarray}
\bar{\Upsilon}(t)&=&R(t)+\bar{B}'(t)\bar{P}(t)\bar{B},\nonumber\\
\bar{\Gamma}(t)&=&B'(t)\bar{P}(t)+\bar{B}'(t)\bar{P}(t)\bar{A}(t).\nonumber
\end{eqnarray}
\end{lemma}
\emph{Proof.} The proof is given in Appendix \ref{a1}.\hfill $\blacksquare$ \\

At the end of this subsection, we give an equivalent solvability condition for Problem 2 by reformulating FBSDEs (\ref{s1}), (\ref{2}) and (\ref{3}) with the denotations below (\ref{jnYC1}).

\begin{lemma}\label{lem2}
If Problem 2 has a solution, then the optimal controller satisfies
\begin{eqnarray}
u(t)&=&-\Upsilon_0^{\dag}(t)\Big[\Gamma_0(t)x(t)+B'(t)\Theta(t)+\bar{B}'(t)\bar{\Theta}(t)\Big]\nonumber\\
&&+G(t)u_1(t),\label{s6}
\end{eqnarray}
where $u_1(t)\in R^{m-m_0}$ is an arbitrary vector such that
\begin{eqnarray}
0&=&C_0(t)x(t)+B_0'(t)\Theta(t)+\bar{B}_0'(t)\bar{\Theta}(t), \label{n2}
\end{eqnarray}
and $(x(t), \Theta(t),\bar{\Theta}(t))$ obey the following FBSDEs:
\begin{eqnarray}
dx(t)&=&\Big[A_0(t)x(t)+D_0(t)\Theta(t)+F_0(t)\bar{\Theta}(t)\nonumber\\
&&+B_{0}(t)u_1(t)\Big]dt+\Big[\bar{A}_0(t)x(t)+\bar{D}_0(t)\Theta(t)\nonumber\\
&&+\bar{F}_0(t)\bar{\Theta}(t)+\bar{B}_0(t)u_1(t)\Big]dw(t),\label{n3}\\
d\Theta(t)&=&-\Big[A'_0(t)\Theta(t)+\bar{A}'_0(t)\bar{\Theta}(t)+C'_0(t) u_1(t)\Big]dt\nonumber\\
&& +\bar{\Theta}(t)dw(t), \label{c8}
\end{eqnarray}
with $x(0)=x_0$ and $\Theta(T)=0.$ Conversely, if FBSDEs (\ref{n3}), (\ref{c8}) and (\ref{n2}) is solvable, then Problem 2 is also solvable.
\end{lemma}
\emph{Proof.} The proof is given in Appendix \ref{a3}. \hfill $\blacksquare$

\subsection{Solution to stochastic optimal control with irregular performance}

We are now in the position to present the main result for the stochastic optimal control with irregular performance.

\begin{theorem}\label{thsto}
Under Assumption \ref{assum2}, Problem 2 is solvable
if there exists a matrix $P_1(T)$ such that the following changed cost function
\begin{eqnarray}
\bar{J}(x_0;u)=J(x_0;u)+E\Big\{x'(T)\Big[H+P_1(T)\Big]x(T)\Big\}, \label{s5}
\end{eqnarray}
is regular and  $P_1(T)x(T)=0$ is achieved with the controller minimizing (\ref{s5}).

\end{theorem}
\emph{Proof.} Based on Lemma \ref{lem2}, Problem 2 is solvable if FBSDEs (\ref{n3}), (\ref{c8}) and (\ref{n2}) is solvable.
Thus, it is sufficient to prove that if there exists a matrix $P_1(T)$ such that (\ref{s5}) is regular and $P_1(T)x(T)=0$, then the FBSDEs (\ref{n3}), (\ref{c8}) and (\ref{n2}) is solvable.
To this end, the proof is divided into two steps. The first step is to show that the regularity of (\ref{s5}) implies
the following condition holds:
\begin{eqnarray}
0&=&C_0(t)+B_{0}'(t)P_1(t)+\bar{B}_{0}'(t)\Big[I-P_1(t)\bar{F}_{0}(t)\Big]^{\dag}\nonumber\\
&&\times P_1(t) \Big[\bar{A}_0(t)+\bar{D}_{0}(t) P_1(t)\Big].\label{pz2}
\end{eqnarray}
The second step is to verify that under the conditions of Assumption 2, (\ref{pz2}) and $P_1(T)x(T)=0$, the following defined $(x(t),\Theta(t),\bar{\Theta}(t))$ solves the FBSDEs (\ref{n3}), (\ref{c8}) and (\ref{n2}):
\begin{eqnarray}
\Theta(t)&=&P_1(t)x(t),\label{p3}\\
\bar{\Theta}(t)&=&P_1(t)\Big[\bar{A}_1(t)x(t)+\bar{B}_1(t)u_1(t)\Big],\label{p4}\\
dx(t)&=&\Big[A_1(t)x(t)+B_1(t)u_1(t)\Big]dt\nonumber\\
&&+\Big[\bar{A}_1(t)x(t)+\bar{B}_{1}(t)u_1(t)\Big]dw(t).\label{p6}
\end{eqnarray}

First of all, we prove that if (\ref{s5}) is regular, then (\ref{pz2}) holds.
In fact, by using (\ref{s9}) in Appendix \ref{a1}, we have
\begin{eqnarray}
&&\Big\{\Gamma_0(t)+B'(t)P_1(t)+\bar{B}'(t)\nonumber\\
&&\times \Big[I-P_1(t)\bar{F}_{0}(t)\Big]^{\dag}P_1(t)\Big[\bar{A}_0(t)+\bar{D}_{0}(t) P_1(t)\Big]\Big\}'\nonumber\\
&=&\Big[\bar{A}'(t)P_1(t)\bar{B}(t)+\Gamma_0'(t)+P_1(t)B(t)\Big]\nonumber\\
&&\times\Big\{I-\Upsilon_0^{\dag}(t)\bar{B}'(t)\Big[I-P_1(t)\bar{F}_{0}(t)\Big]^{\dag} P_1(t)\bar{B}(t)\Big\}\nonumber\\
&=&\bar{\Gamma}'(t)\Big\{I-\Upsilon_0^{\dag}(t)\bar{B}'(t)\Big[I-P_1(t)\bar{F}_{0}(t)\Big]^{\dag} P_1(t)\bar{B}(t)\Big\}.\nonumber
\end{eqnarray}
Together with the denotations below (\ref{jnYC1}), it is further obtained that
\begin{eqnarray}
&&\Upsilon_{T_0}'(t)\Big\{C_0(t)+B_{0}'(t)P_1(t)+\bar{B}_{0}'(t)\Big[I-P_1(t)\bar{F}_{0}(t)\Big]^{\dag}\nonumber\\
&&\times P_1(t) \Big[\bar{A}_0(t)+\bar{D}_{0}(t) P_1(t)\Big]\Big\}\nonumber\\
&=&\Big[I-\Upsilon_0(t)\Upsilon_0^{\dag}(t)\Big]\Big\{I-\Upsilon_0^{\dag}(t)\bar{B}'(t) \Big[I-P_1(t)\bar{F}_{0}(t)\Big]^{\dag} \nonumber\\
&&\times P_1(t)\bar{B}(t)\Big\}'\bar{\Gamma}(t)\nonumber\\
&=&\Big\{I-\Upsilon_0(t)\bar{\Upsilon}^{\dag}(t)-\bar{B}'(t)P_1(t)\Big[I-\bar{F}_{0}(t)P_1(t)\Big]^{\dag}\nonumber\\
&&\times \bar{B}(t)\Upsilon_0^{\dag}(t)\Big\}\bar{\Gamma}(t)\nonumber\\
&=&\Big[I-\bar{\Upsilon}(t)\bar{\Upsilon}^{\dag}(t)\Big]\bar{\Gamma}(t)\nonumber\\
&=&0,\label{s8}
\end{eqnarray}
where the last equality follows from the regularity of (\ref{s5}). In view of the fact that $\Upsilon_{T_0}(t)$ has full row rank,
we obtain (\ref{pz2}) holds.

Next, we prove that $(x(t),\Theta(t),\bar{\Theta}(t))$ defined by (\ref{p3})-(\ref{p6}) solves the FBSDEs (\ref{n3}), (\ref{c8}) and (\ref{n2}). In fact, by taking It\^{o}'s formula to $P_1(t)x(t)$, we have
\begin{eqnarray}
&&d\Big[P_1(t)x(t)\Big]\nonumber\\
&=&\dot{P}_1(t)x(t)dt+P_1(t)\Big[A_1(t)x(t)+B_1(t)u_1(t)\Big]dt\nonumber\\
&&+P_1(t)\Big[\bar{A}_1(t)x(t)+\bar{B}_{1}(t)u_1(t)\Big]dw(t)\nonumber\\
&=&-A_0'(t)P_1(t)x(t)dt-\bar{A}_0'(t)\Big[I-P_1(t)\bar{F}_{0}(t)\Big]^{\dag}\nonumber\\
&&\times \Big\{P_1(t)\Big[\bar{A}_0(t)+\bar{D}_0(t)P_1(t)\Big]x(t)\nonumber\\
&&+P_1(t)\bar{B}_0(t)u_1(t)\Big\}dt\nonumber\\
&&+\Big\{P_1(t)B_0(t)+\Big[\bar{A}_0'(t)+P_1(t)F_0(t)\Big]\nonumber\\
&&\times \Big[I-P_1(t)\bar{F}_{0}(t)\Big]^{\dag}P_1(t)\bar{B}_0(t)\Big\}u_1(t)dt\nonumber\\
&&+P_1(t)\Big[\bar{A}_1(t)x(t)+\bar{B}_1(t)u_1(t)\Big]dw(t),\nonumber
\end{eqnarray}
where the equation (\ref{p1}) of $P_1(t)$ has been used in the derivation of the last equality.
Combining with (\ref{pz2}), the above equation is further formulated as
\begin{eqnarray}
&&d\Big[P_1(t)x(t)\Big]\nonumber\\
&=&-A_0'(t)P_1(t)x(t)dt-\bar{A}_0'(t)\Big[I-P_1(t)\bar{F}_{0}(t)\Big]^{\dag}\nonumber\\
&&\times \Big\{P_1(t)\Big[\bar{A}_0(t)+\bar{D}_0(t)P_1(t)\Big]x(t)+P_1(t)\bar{B}_0(t)\nonumber\\
&&\times u_1(t)\Big\}dt-C_0'(t)u_1(t)dt\nonumber\\
&&+P_1(t)\Big[A_1(t)x(t)+B_1(t)u_1(t)\Big]dw(t).\label{p7}
\end{eqnarray}
By taking again (\ref{s9}) in Appendix \ref{a1}, we obtain that
\begin{eqnarray}
&&L_1'(t)\Big[I-P_1(t)\bar{F}_{0}(t)\Big]^{\dag}P_1(t)L_2(t)\nonumber\\
&=&L_1'(t)\Big\{I+P_1(t)\bar{F}_{0}(t) \Big[I-P_1(t)\bar{F}_{0}(t)\Big]^{\dag}\Big\}\nonumber\\
&&\times P_1(t)L_2(t).\nonumber
\end{eqnarray}
Accordingly, (\ref{p7}) can be reformulated as
\begin{eqnarray}
&&d\Big[P_1(t)x(t)\Big]\nonumber\\
&=&-A_0'(t)P_1(t)x(t)dt-\bar{A}_0'(t)\Big\{I+P_1(t)\bar{F}_{0}(t)\nonumber\\
&&\times \Big[I-P_1(t)\bar{F}_{0}(t)\Big]^{\dag}\Big\} \Big\{P_1(t)\Big[\bar{A}_0(t)\nonumber\\
&&+\bar{D}_0(t)P_1(t)\Big]x(t)+P_1(t)\bar{B}_0(t) u_1(t)\Big\}dt\nonumber\\
&&-C_0'(t)u_1(t)dt+P_1(t)\Big[\bar{A}_1(t)x(t)\nonumber\\
&&+\bar{B}_1(t)u_1(t)\Big]dw(t)\nonumber\\
&=&-A_0'(t)P_1(t)x(t)dt-\bar{A}_0'(t)P_1(t)\Big[\bar{A}_1(t)x(t)\nonumber\\
&&+\bar{B}_1(t)u_1(t)\Big]dt-C_0'(t)u_1(t)dt\nonumber\\
&&+P_1(t)\Big[\bar{A}_1(t)x(t)+\bar{B}_1(t)u_1(t)\Big]dw(t).\label{p5}
\end{eqnarray}
In addition, by using (\ref{pz1}) and (\ref{pz2}), we have
\begin{eqnarray}
&&C_0(t)x(t)+B_0'(t)P_1(t)x(t)\nonumber\\
&&+\bar{B}_0'(t)P_1(t)\Big[\bar{A}_1(t)x(t)+\bar{B}_1(t)u_1(t)\Big]\nonumber\\
&=&\Big\{C_0(t)+B_0'(t)P_1(t)+\bar{B}_0'(t)\Big[I-P_1(t)\bar{F}_{0}(t)\Big]^{\dag}\nonumber\\
&&\times P_1(t)\Big[\bar{A}_0(t)+\bar{D}_0(t)P_1(t)\Big]\Big\}x(t)\nonumber\\
&&+\bar{B}_0'(t)\Big[I-P_1(t)\bar{F}_{0}(t)\Big]^{\dag}P_1(t)\bar{B}_0(t)u_1(t)\nonumber\\
&=&0.\label{n4}
\end{eqnarray}
Similarly, we can rewrite (\ref{p6}) as follows:
\begin{eqnarray}
dx(t)&=&\Big\{A_0(t)x(t)+D_0(t)P_1(t)x(t)+F_0(t)P_1(t)\nonumber\\
&&\times \Big[\bar{A}_1(t)x(t)+\bar{B}_1(t)u_1(t)\Big]+B_{0}(t)u_1(t)\Big\}dt\nonumber\\
&&+\Big\{\bar{A}_0(t)x(t)+\bar{D}_0(t)P_1(t)x(t)+\bar{F}_0(t)P_1(t)\nonumber\\
&&\times \Big[\bar{A}_1(t)x(t)+\bar{B}_1(t)u_1(t)\Big]+\bar{B}_0(t)u_1(t)\Big\}dw(t).\nonumber\\\label{p9}
\end{eqnarray}
By making comparison between (\ref{p9}), (\ref{p5}), (\ref{n4}) and (\ref{n3}), (\ref{c8}) and (\ref{n2}), and using the existence of $u_1(t)$ such that $P_1(T)x(T)=0$, it follows that
(\ref{p3})-(\ref{p6}) solves the FBSDEs (\ref{n3}), (\ref{c8}) and (\ref{n2}).

Based on Lemma \ref{lem2}, Problem 2 is solvable. The proof is now completed. \hfill $\blacksquare$\\

%

We now present the optimal controller for the stochastic optimal control with irregular performance.
\begin{corollary}
Under the Assumption \ref{assum2}, if there exists a matrix $P_1(T)$ such that the cost function (\ref{s5}) is regular and $P_1(T)x(T)=0$ can be achieved with the controller minimizing (\ref{s5}),
then the optimal controller $u(t)$ is given by
\begin{eqnarray}
u(t)&=&-\bar{\Upsilon}^{\dag}(t)\bar{\Gamma}(t)x(t)+[I-\bar{\Upsilon}^{\dag}(t)\bar{\Upsilon}(t)]z(t),\label{r2}
\end{eqnarray}
where $z(t)\in R^m$ is an arbitrary vector with compatible dimension such that $P_1(T)x(T)=0$.
The optimal cost is given by
\begin{eqnarray}
J^*(x_0;u)&=&x_0'\bar{P}(t_0)x_0.\label{p8}
\end{eqnarray}
\end{corollary}
\emph{Proof.} By solving the regular optimal control problem of minimizing (\ref{s5}) subject to (\ref{s1}), we have
the optimal control (\ref{r2}) directly. 
The proof is now completed.
\hfill $\blacksquare$

\section{Conclusions}

In this paper, we have investigated the essential problem of the irregular LQ control.
It was shown that the difference between the irregular LQ control and the standard (regular) one is that the irregular controller needs to do two things at the same time (minimizing the LQ cost and achieving the state terminal condition).  As application, we have presented a sufficient condition solution to the irregular LQ control for stochastic systems with multiplicative noise.


\appendix

\subsection{Proof of Lemma \ref{leminv}}\label{a2}
\begin{enumerate}
  \item
  By using (\ref{p11}), we have
  \begin{eqnarray}
  &&L_1'(t)P_1(t)\Big[I-\bar{F}_0(t)P_1(t)\Big]^{\dag}L_2(t)\nonumber\\
  &=&L_1'(t)\Big[I-P_1(t)\bar{F}_{0}(t)\Big]^{\dag}\Big[I-P_1(t)\bar{F}_{0}(t)\Big]P_1(t)\nonumber\\
  &&\times  \Big[I-\bar{F}_0(t)P_1(t)\Big]^{\dag}L_2(t)\nonumber\\
  &=&L_1'(t)\Big[I-P_1(t)\bar{F}_{0}(t)\Big]^{\dag}P_1(t)\nonumber\\
  &&\times \Big[I-\bar{F}_{0}(t)P_1(t)\Big]\Big[I-\bar{F}_0(t)P_1(t)\Big]^{\dag}L_2(t)\nonumber\\
  &=&L_1'(t)\Big[I-P_1(t)\bar{F}_{0}(t)\Big]^{\dag}P_1(t)L_2(t).\nonumber
  \end{eqnarray}

  \item  By using again (\ref{p11}), it is obtained
  \begin{eqnarray}
&&\Big[I-\bar{F}_0(t)P_1(t)\Big]^{\dag}L(t)\nonumber\\
&=&L(t)+\bar{F}_0(t)P_1(t)\Big[I-\bar{F}_0(t)P_1(t)\Big]^{\dag}L(t)\nonumber\\
&=&L(t)+\bar{F}_0(t)\Big[I-P_1(t)\bar{F}_{0}(t)\Big]^{\dag}\Big[I-P_1(t)\bar{F}_{0}(t)\Big]\nonumber\\
&&\times P_1(t) \Big[I-\bar{F}_0(t)P_1(t)\Big]^{\dag}L(t)\nonumber\\
&=&L(t)+\bar{F}_0(t)\Big[I-P_1(t)\bar{F}_{0}(t)\Big]^{\dag}P_1(t)\nonumber\\
&&\times\Big[I-\bar{F}_{0}(t)P_1(t)\Big] \Big[I-\bar{F}_0(t)P_1(t)\Big]^{\dag}L(t)\nonumber\\
&=&L(t)+\bar{F}_0(t)\Big[I-P_1(t)\bar{F}_{0}(t)\Big]^{\dag}P_1(t)L(t).\nonumber
  \end{eqnarray}
  \end{enumerate}

 The proof is now completed.

\subsection{Proof of Lemma \ref{lem3}}\label{a1}

First, we make some algebraic calculations to Riccati equation (\ref{p1}). Based on Lemma \ref{leminv}, it can be derived that
\begin{eqnarray}
L_1'(t)P_1(t)L_2(t)&=&L_1'(t)\Big[I+P_1(t)\bar{B}(t)\Upsilon_0^{\dag}(t)\bar{B}'(t)\Big]\nonumber\\
&&\hspace{-26mm}\times  \Big[I+P_1(t)\bar{B}(t)\Upsilon_0^{\dag}(t)\bar{B}'(t)\Big]^{\dag}P_1(t)L_2(t),\label{s9}
\end{eqnarray}
where $L_1(t), L_2(t)$ may be $B(t), \bar{A}(t)$, or $\bar{B}(t)$.
This implies that
\begin{eqnarray}
&&\bar{A}'(t)\Big[I+P_1(t)\bar{B}(t)\Upsilon_0^{\dag}(t)\bar{B}'(t)\Big]^{\dag} P_1(t)\bar{B}(t)\Upsilon_0^{\dag}(t)\nonumber\\
&&\times \Gamma_{0}(t)\nonumber\\
&=&\bar{A}'(t)P_1(t)\bar{B}(t)\Upsilon_0^{\dag}(t)\Gamma_{0}(t)\nonumber\\
&&-\bar{A}'(t)P_1(t)\bar{B}(t)\Upsilon_0^{\dag}(t)\bar{B}'(t)\nonumber\\
&&\times \Big[I+P_1(t)\bar{B}(t)\Upsilon_0^{\dag}(t)\bar{B}'(t)\Big]^{\dag}\nonumber\\
&&\times  P_1(t)\bar{B}(t)\Upsilon_0^{\dag}(t)\Gamma_{0}(t)\nonumber\\
&=&\bar{A}'(t)P_1(t)\bar{B}(t)\Big[\Upsilon_0(t)+\bar{B}'(t)P_1(t)\bar{B}(t)\Big]^{\dag}\Gamma_{0}(t),\label{a4}
\end{eqnarray}
where (\ref{p2}) has been used in the derivation of the last equality. By taking similar procedures to the above equation,
it is obtained that
\begin{eqnarray}
&&\Big[\bar{A}'_0(t)+P_1(t)F_0(t)\Big]\Big[I-P_1(t)\bar{F}_{0}(t)\Big]^{\dag}\nonumber\\
&&\times P_1(t)\Big[\bar{A}_0(t)+\bar{D}_0(t)P_1(t)\Big]\nonumber\\
&=&\Big[\bar{A}'(t)-\Gamma_{0}'(t)\Upsilon_0^{\dag}(t)\bar{B}'(t)-P_1(t)B(t)\Upsilon_0^{\dag}(t)\bar{B}'(t)\Big]\nonumber\\
&&\times \Big[I+P_1(t)\bar{B}(t)\Upsilon_0^{\dag}(t)\bar{B}'(t)\Big]^{\dag} P_1(t)\nonumber\\
&&\times\Big[\bar{A}(t)-\bar{B}(t)\Upsilon_0^{\dag}(t)\Gamma_{0}(t)-\bar{B}(t)\Upsilon_0^{\dag}(t)B'(t)P_1(t)\Big]\nonumber\\
&=&\bar{A}'(t)\Big[I+P_1(t)\bar{B}(t)\Upsilon_0^{\dag}(t)\bar{B}'(t)\Big]^{\dag} P_1(t)\bar{A}(t)\nonumber\\
&&-\bar{A}'(t)\Big[I+P_1(t)\bar{B}(t)\Upsilon_0^{\dag}(t)\bar{B}'(t)\Big]^{\dag} P_1(t)\bar{B}(t)\Upsilon_0^{\dag}(t)\nonumber\\
&&\times \Gamma_{0}(t)\nonumber\\
&&-\bar{A}'(t)\Big[I+P_1(t)\bar{B}(t)\Upsilon_0^{\dag}(t)\bar{B}'(t)\Big]^{\dag} P_1(t)\bar{B}(t)\Upsilon_0^{\dag}(t)\nonumber\\
&&\times B'(t)P_1(t)\nonumber\\
&&-\Gamma_{0}'(t)\Upsilon_0^{\dag}(t)\bar{B}'(t)\Big[I+P_1(t)\bar{B}(t)\Upsilon_0^{\dag}(t)\bar{B}'(t)\Big]^{\dag}\nonumber\\
&&\times  P_1(t)\bar{A}(t)\nonumber\\
&&+\Gamma_{0}'(t)\Upsilon_0^{\dag}(t)\bar{B}'(t)\Big[I+P_1(t)\bar{B}(t)\Upsilon_0^{\dag}(t)\bar{B}'(t)\Big]^{\dag} P_1(t)\nonumber\\
&&\times \bar{B}(t)\Upsilon_0^{\dag}(t)\Gamma_{0}(t)\nonumber\\
&&+\Gamma_{0}'(t)\Upsilon_0^{\dag}(t)\bar{B}'(t)\Big[I+P_1(t)\bar{B}(t)\Upsilon_0^{\dag}(t)\bar{B}'(t)\Big]^{\dag} P_1(t)\nonumber\\
&&\times \bar{B}(t)\Upsilon_0^{\dag}(t)B'(t)P_1(t)\nonumber\\
&&-P_1(t)B(t)\Upsilon_0^{\dag}(t)\bar{B}'(t)\Big[I+P_1(t)\bar{B}(t)\Upsilon_0^{\dag}(t)\bar{B}'(t)\Big]^{\dag} \nonumber\\
&&\times P_1(t)\bar{A}(t)\nonumber\\
&&+P_1(t)B(t)\Upsilon_0^{\dag}(t)\bar{B}'(t)\Big[I+P_1(t)\bar{B}(t)\Upsilon_0^{\dag}(t)\bar{B}'(t)\Big]^{\dag}\nonumber\\
&&\times P_1(t) \bar{B}(t)\Upsilon_0^{\dag}(t)\Gamma_{0}(t)\nonumber\\
&&+P_1(t)B(t)\Upsilon_0^{\dag}(t)\bar{B}'(t)\Big[I+P_1(t)\bar{B}(t)\Upsilon_0^{\dag}(t)\bar{B}'(t)\Big]^{\dag} \nonumber\\
&&\times P_1(t) \bar{B}(t)\Upsilon_0^{\dag}(t)B'(t)P_1(t)\nonumber\\
&=&\bar{A}'(t)P_1(t)\bar{A}(t)\nonumber\\
&&-\bar{A}'(t)P_1(t)\bar{B}(t)\Big[\Upsilon_0(t)+\bar{B}'(t)P_1(t)\bar{B}(t)\Big]^{\dag}\nonumber\\
&&\times  \bar{B}'(t)P_1(t)\bar{A}(t)\nonumber\\
&&-\bar{A}'(t)P_1(t)\bar{B}(t)\Big[\Upsilon_0(t)+\bar{B}'(t)P_1(t)\bar{B}(t)\Big]^{\dag} \Gamma_{0}(t)\nonumber\\
&&-\bar{A}'(t)P_1(t)\bar{B}(t)\Big[\Upsilon_0(t)+\bar{B}'(t)P_1(t)\bar{B}(t)\Big]^{\dag} \nonumber\\
&&\times B'(t)P_1(t)\nonumber\\
&&-\Gamma_{0}'(t)\Big[\Upsilon_0(t)+\bar{B}'(t)P_1(t)\bar{B}(t)\Big]^{\dag}\nonumber\\
&&\times \bar{B}'(t)P_1(t)\bar{A}(t)\nonumber\\
&&-\Gamma_{0}'(t)\Big[\Upsilon_0(t)+\bar{B}'(t)P_1(t)\bar{B}(t)\Big]^{\dag}\Gamma_{0}(t)\nonumber\\
&&+\Gamma_{0}'(t)\Upsilon_0^{\dag}(t)\Gamma_{0}(t)\nonumber\\
&&-\Gamma_{0}'(t)\Big[\Upsilon_0(t)+\bar{B}'(t)P_1(t)\bar{B}(t)\Big]^{\dag}B'(t)P_1(t)\nonumber\\
&&+\Gamma_{0}'(t)\Upsilon_0^{\dag}(t)B'(t)P_1(t)\nonumber\\
&&-P_1(t)B(t)\Big[\Upsilon_0(t)+\bar{B}'(t)P_1(t)\bar{B}(t)\Big]^{\dag}\nonumber\\
&&\times \bar{B}'(t) P_1(t)\bar{A}(t)\nonumber\\
&&-P_1(t)B(t)\Big[\Upsilon_0(t)+\bar{B}'(t)P_1(t)\bar{B}(t)\Big]^{\dag}\Gamma_{0}(t)\nonumber\\
&&+P_1(t)B(t)\Upsilon_0^{\dag}(t)\Gamma_{0}(t)\nonumber\\
&&-P_1(t)B(t)\Big[\Upsilon_0(t)+\bar{B}'(t)P_1(t)\bar{B}(t)\Big]^{\dag}\nonumber\\
&&\times B'(t)P_1(t)\nonumber\\
&&+P_1(t)B(t)\Upsilon_0^{\dag}(t)B'(t)P_1(t)\nonumber\\
&=&\bar{A}'(t)P_1(t)\bar{A}(t)\nonumber\\
&&-\Big[\bar{A}'(t)P_1(t)\bar{B}(t)+P_1(t)B(t)+\Gamma_{0}'(t)\Big]\nonumber\\
&&\times \Big[\Upsilon_0(t)+\bar{B}'(t)P_1(t)\bar{B}(t)\Big]^{\dag}\nonumber\\
&&\times \Big[\bar{B}'(t)P_1(t)\bar{A}(t)+B'(t)P_1(t)+\Gamma_{0}(t)\Big]\nonumber\\
&&+\Gamma_{0}'(t)\Upsilon_0^{\dag}(t)\Gamma_{0}(t)\nonumber\\
&&+\Gamma_{0}'(t)\Upsilon_0^{\dag}(t)B'(t)P_1(t)\nonumber\\
&&+P_1(t)B(t)\Upsilon_0^{\dag}(t)\Gamma_{0}(t)\nonumber\\
&&+P_1(t)B(t)\Upsilon_0^{\dag}(t)B'(t)P_1(t).\nonumber
\end{eqnarray}
Together with  (\ref{r1}) and (\ref{p1}), it follows that
\begin{eqnarray}
0&=&\dot{P}_1(t)+\dot{P}(t)+\Big[P_1(t)+P(t)\Big]A(t)\nonumber\\
&&+A'(t)\Big[P_1(t)+P(t)\Big]\nonumber\\
&&-P_1(t)B(t)\Upsilon_0^{\dag}(t)\Gamma_{0}(t)\nonumber\\
&&-\Gamma_{0}'(t)\Upsilon_0^{\dag}(t)B'(t)P_1(t)\nonumber\\
&&-P_1(t)B(t)\Upsilon_0^{\dag}(t)B'(t)P_1(t)\nonumber\\
&&+\bar{A}'(t)\Big[P_1(t)+P(t)\Big]\bar{A}(t)\nonumber\\
&&-\Big[\bar{A}'(t)P_1(t)\bar{B}(t)+P_1(t)B(t)+\Gamma_{0}'(t)\Big]\nonumber\\
&&\times \Big[\Upsilon_0(t)+\bar{B}'(t)P_1(t)\bar{B}(t)\Big]^{\dag}\nonumber\\
&&\times \Big[\bar{B}'(t)P_1(t)\bar{A}(t)+B'(t)P_1(t)+\Gamma_{0}(t)\Big]\nonumber\\
&&+\Gamma_{0}'(t)\Upsilon_0^{\dag}(t)\Gamma_{0}(t)\nonumber\\
&&+\Gamma_{0}'(t)\Upsilon_0^{\dag}(t)B'(t)P_1(t)\nonumber\\
&&+P_1(t)B(t)\Upsilon_0^{\dag}(t)\Gamma_{0}(t)\nonumber\\
&&+P_1(t)B(t)\Upsilon_0^{\dag}(t)B'(t)P_1(t)\nonumber\\
&&-\Gamma_0'(t)\Upsilon_0^{\dag}(t)\Gamma_0(t)\nonumber\\
&=&\dot{P}_1(t)+\dot{P}(t)+\Big[P_1(t)+P(t)\Big]A(t)\nonumber\\
&&+A'(t)\Big[P_1(t)+P(t)\Big]\nonumber\\
&&+\bar{A}'(t)\Big[P_1(t)+P(t)\Big]\bar{A}(t)\nonumber\\
&&-\Big[\bar{A}'(t)P_1(t)\bar{B}(t)+P_1(t)B(t)+\Gamma_{0}'(t)\Big]\nonumber\\
&&\times \Big[R(t)+\bar{B}'(t)(P_1(t)+P_1(t))\bar{B}(t)]^{\dag}\nonumber\\
&&\times \Big[\bar{B}'(t)P_1(t)\bar{A}(t)+B'(t)P_1(t)+\Gamma_{0}(t)\Big].\nonumber
\end{eqnarray}
This is exactly (\ref{s7}). The proof is now completed.

\subsection{Proof of Lemma \ref{lem2}}\label{a3}

The key is to reformulate FBSDEs (\ref{s1}), (\ref{2}) and (\ref{3}) as FBSDEs (\ref{n3}), (\ref{c8}) and (\ref{n2}). Without loss of generality, we assume that
\begin{eqnarray}
\Theta(t)= p(t)-P(t)x(t),\label{21}
\end{eqnarray}
where $P(t)$ obeys the Riccati equation (\ref{r1}). It is obvious that $P(T)=H$ and $\Theta(T)=p(T)-P(T)x(T)=0$.
We also assume without loss of generality that
\begin{eqnarray}
d\Theta(t)=\hat{\Theta}(t)dt+\bar{\Theta}(t)dw(t), \label{zj1}
\end{eqnarray}
where $\hat{\Theta}(t)$ and $\bar{\Theta}(t)$ are to be determined.  Applying It\^{o}'s formula to (\ref{21}) yields
\begin{eqnarray}
dp(t)&=&\dot{P}(t)x(t)dt+P(t)\Big[A(t)x(t)
+B(t)u(t)\Big]dt\nonumber\\
&&+P(t)\Big[\bar{A}(t)x(t)+\bar{B}(t)u(t)\Big]dw(t)+\hat{\Theta}(t) dt\nonumber\\
&&+\bar{\Theta}(t)dw(t).\label{zz4}
\end{eqnarray}
Using (\ref{21}), we rewrite (\ref{2}) as
\begin{eqnarray}
dp(t)&=&-\Big[A'(t)P(t)x(t)+A'(t)\Theta(t)
+\bar{A}'(t)q(t)\nonumber\\
&&+Q(t)x(t)\Big]dt+q(t)dw(t).\label{11}
\end{eqnarray}
With a comparison of (\ref{zz4}) and (\ref{11}), it follows that
\begin{eqnarray}
q(t)&=& P(t)\Big[\bar{A}(t)x(t)+\bar{B}(t)u(t)\Big]+\bar{\Theta}(t),\label{zz1}\\
0&=&\dot{P}(t)x(t)+P(t)A(t)x(t)+P(t)B(t)u(t)\nonumber\\
&&+\hat{\Theta}(t)
+A'(t)P(t)x(t)+A'(t)\Theta(t)\nonumber\\
&&+\bar{A}'(t)q(t)+Q(t)x(t). \label{r6}
\end{eqnarray}
Using (\ref{zz1}) and (\ref{21}), (\ref{3}) becomes
\begin{eqnarray}
0&=&R(t)u(t)+B'(t)P(t)x(t)+B'(t)\Theta(t)\nonumber\\
&&+\bar{B}'(t)P(t)\bar{A}(t)x(t)+\bar{B}'(t)P(t)\bar{B}(t)u(t)\nonumber\\
&&+\bar{B}'(t)\bar{\Theta}(t)\nonumber\\
&=&\Upsilon_0(t)u(t)
+\Gamma_0(t)x(t)+B'(t)\Theta(t)\nonumber\\
&&+\bar{B}'(t)\bar{\Theta}(t),\label{25}
\end{eqnarray}
where $\Upsilon_0(t)$ and  $\Gamma_0(t)$ are respectively as in (\ref{Intr3}) and (\ref{Intr4}).

Thus, (\ref{25}) can be equivalently written as
\begin{eqnarray}
u(t)&=&-\Upsilon_0^{\dag}(t)\Big[\Gamma_0(t)x(t)+B'(t)\Theta(t)+\bar{B}'(t)\bar{\Theta}(t)\Big]\nonumber\\
&&+\Big[I-\Upsilon_0^{\dag}(t)\Upsilon_0(t)\Big]z(t),\label{n1}
\end{eqnarray}
where $z(t)$ is a vector with compatible dimension such that the following equality holds
\begin{eqnarray}
0&=&\Big[I-\Upsilon_0(t)\Upsilon_0^{\dag}(t)\Big]\Big[\Gamma_0(t)x(t)+B'(t)\Theta(t)\nonumber\\
&&+\bar{B}'(t)\bar{\Theta}(t)\Big]. \label{nz1}
\end{eqnarray}
Denote \begin{eqnarray}
T_0(t)\Big[I-\Upsilon_0^{\dag}(t)\Upsilon_0(t)\Big]z(t)=\left[
                              \begin{array}{c}
                                 0  \\
                                 u_1(t)\\
                              \end{array}
                            \right], \label{jn1}
\end{eqnarray}
where $u_1(t)=\Upsilon_{T_0}(t) z(t)\in R^{m-m_0}$. Note that $\Upsilon_{T_0}(t)\in R^{[m-m_0]\times m}$ is full row rank, it yields that $u_1(t)$ can be arbitrary due to
the arbitrariness of $z(t)$. Moreover, the optimal controller (\ref{s6}) follows from (\ref{n1}) by using the denotation of $G(t)$.

Now we rewrite (\ref{nz1}) as (\ref{n2}). First, it is noted that
\begin{eqnarray}
&&I-\Upsilon_0(t)\Upsilon_0^{\dag}(t)\nonumber\\
&=&\Big[I-\Upsilon_0(t)\Upsilon_0^{\dag}(t)\Big]\Big[I-\Upsilon_0(t)\Upsilon_0^{\dag}(t)\Big]\nonumber \\
&=&\Big[I-\Upsilon_0(t)\Upsilon_0^{\dag}(t)\Big]T'_0(t)\Big[T^{-1}_0(t)\Big]'\Big[I-\Upsilon_0(t)\Upsilon_0^{\dag}(t)\Big]\nonumber \\
&=&\left[
\begin{array}{cc}
                                 0 &
                                 \Upsilon'_{T_0}(t)\\
                              \end{array}
                            \right]\Big[T^{-1}_0(t)\Big]'\Big[I-\Upsilon_0(t)\Upsilon_0^{\dag}(t)\Big],
\end{eqnarray}
where (\ref{jnYC1}) has been used in the derivation of the last equality. By using the denotations below (\ref{jnYC1}), (\ref{nz1}) can be written as
\begin{eqnarray}
0&=&\Upsilon'_{T_0}(t)\Big[C_0(t)x(t)+B_0'(t)\Theta(t)+\bar{B}_0'(t)\bar{\Theta}(t)\Big].\label{nYc2}
\end{eqnarray}
Note that $\Upsilon'_{T_0}(t)$ is full column rank, (\ref{nYc2}) is rewritten as (\ref{n2}) directly.
By substituting (\ref{n1}) and (\ref{zz1}) into (\ref{r6}) and using (\ref{r1}), it yields that
\begin{eqnarray}
0&=&\dot{P}(t)x(t)+P(t)A(t)x(t)+\hat{\Theta}(t)
+A'(t)P(t)x(t)\nonumber\\
&&+A'(t)\Theta(t)+\bar{A}'(t)P(t)\bar{A}(t)x(t)+\bar{A}'(t)\bar{\Theta}(t)\nonumber\\
&&+Q(t)x(t)-\Gamma_0'(t)\Upsilon_0^{\dag}(t)\Big[\Gamma_0(t)x(t)+B'(t)\Theta(t)\nonumber\\
&&+\bar{B}'(t)\bar{\Theta}(t)\Big]+\Gamma_0'(t)\Big[I-\Upsilon_0^{\dag}(t)\Upsilon_0(t)\Big]z(t)\nonumber\\
&=&\hat{\Theta}(t)+\Big[A'(t)-\Gamma_0'(t)\Upsilon_0^{\dag}(t)B'(t)\Big]\Theta(t)\nonumber\\
&&+\Big[\bar{A}'(t)-\Gamma_0'(t)\Upsilon_0^{\dag}(t)\bar{B}'(t)\Big]\bar{\Theta}(t)\nonumber\\
&&+\Gamma_0'(t)\Big[I-\Upsilon_0^{\dag}(t)\Upsilon_0(t)\Big]z(t).\label{ZJ1}
\end{eqnarray}

In view of the fact that $[I-\Upsilon_0^{\dag}(t)\Upsilon_0(t)]^2=I-\Upsilon_0^{\dag}(t)\Upsilon_0(t),$ it is obtained that
\begin{eqnarray}
&&\Gamma_0'(t)\Big[I-\Upsilon_0^{\dag}(t)\Upsilon_0(t)\Big]z(t)\nonumber\\
&=&\Gamma_0'(t)\Big[I-\Upsilon_0^{\dag}(t)\Upsilon_0(t)\Big]T_0^{-1}(t)T_0(t)\nonumber\\
&&\times\Big[I-\Upsilon_0^{\dag}(t)\Upsilon_0(t)\Big]z(t)\nonumber\\
&=&\Gamma_0'(t)\Big[I-\Upsilon_0^{\dag}(t)\Upsilon_0(t)\Big]T_0^{-1}(t)\left[
                              \begin{array}{c}
                                 0  \\
                                 u_1(t)\\
                              \end{array}
                            \right]\nonumber\\
&=&C_0'(t)u_1(t).\label{YCJ1}
\end{eqnarray}
Thus, we have from ({\ref{ZJ1}) that
\begin{eqnarray}
\hat{\Theta}(t)=-\Big[A'_0(t)\Theta(t)+\bar{A}'_0(t)\bar{\Theta}(t)+C'_0(t) u_1(t)\Big],
\end{eqnarray}
this implies that the dynamic of $\Theta(t)$ is as (\ref{c8}) using (\ref{zj1}).

By substituting (\ref{n1}) into (\ref{s1}),
one has the dynamic (\ref{n3}) of the state.

Based on the above derivations, FBSDEs (\ref{s1}), (\ref{2}) and (\ref{3}) has been equivalently rewritten as FBSDEs (\ref{n3}), (\ref{c8}) and (\ref{n2}). Combining with Lemma \ref{lem1},
the result is derived. The proof is now completed.


\ifCLASSOPTIONcaptionsoff
  \newpage
\fi


\begin{thebibliography}{99}


\bibitem{anderson} B. D. O. Anderson, J. B. Moore. \emph{Optimal control: linear quadratic methods}.
Englewood Cliffs, NJ: Prentice Hall, 1990.

\bibitem{Bell} D. J. Bell, Singular problems in optimal control-a survey, \emph{International Journal of Control}, 21(2): 319-331, 1975.


\bibitem{bellman} R. Bellman, The theory of dynamic programming.
\emph{Bulletin of the American Mathematical Society}, 60(6): 503-516, 1954.


\bibitem{Bellman} R. Bellman, I. Glicksberg, O. Gross, \emph{Some aspects of the mathematical theory of control processes}, Rand Corporation, R-313, 1958.

\bibitem{Bonnans} J. F. Bonnans, F. J. Silva, First and second order necessary conditions for
stochastic optimal control problems. Applied Mathematics \& Optimization, 2012, 65: 403-439.


\bibitem{chenhanfu} H.-F. Chen, Unified controls applicable to general case under quadratic index, Acta Mathematicae Applicatae Sinica, 5(1): 45-52, 1982.

\bibitem{chenlizhou} S. Chen, X. Li, X. Zhou, Stochastic linear quadratic regulators with indefinite control weight costs, \emph{SIAM Journal on Control and Optimization}, 36(5): 1685-1702, 1998.

\bibitem{clementsanderson} D. Clements, B. Anderson, \emph{Singular optimal control: The linear-quadratic problem}, Springer-Verlag, New York, 1978.


\bibitem{Gabasov} R. Gabasov, F. M. Kirillova, High order necessary conditions for optimality. SIAM J. Control, 1972, 10: 127-168.

\bibitem{Gurman} V. Gurman, The method of multiple maxima and optimization problems for space maneuvers. Proc. Second Readings of K. E. Tsiolkovskii, Moscow, 1968, 39-51.



\bibitem{Hoehener} D. Hoehener, Variational approach to second-order optimality conditions for control
problems with pure state constraints. SIAM Journal of Control and Optimization, 2012, 50: 1139-1173.



\bibitem{Ho} Y. Ho, Linear Stochastic Singular Control Problems, \emph{Journal of Optimization Theory and Application}, 9(1): 24-31, 1972.


\bibitem{Hsia} T. Hsia, On the existence and synthesis of optimal singular control with
quadratic performance index, \emph{IEEE Trans.\ Autom.\ Control}, 12(6): 778-779, 1967.


\bibitem{Kalman} R. E. Kalman, \emph{Contributions to the theory of optimal control}, Bol. Soc., Mat. Mexicana, 5: 102-119, 1960.

\bibitem{Kliger} I. Kliger, Discussion on the stability of the singular trajectory with respect to ``Bang-Bang" control, \emph{IEEE Trans.\ Autom.\ Control}, 9(4): 583-585, 1964.


\bibitem{Krener}  A. J. Krener, The high order maximal principle and its application to singular extremals. SIAM Journal of Control and Optimization, 1977, 15: 256-293.



\bibitem{Letov} A. M. Letov, The analytical design of control systems, \emph{Automat. Remote Control}, 22: 363-372, 1961.



\bibitem{lewis} F. L. Lewis, D. L. Vrabie, V. L. Syrmos \emph{Optimal control}.
John Wiley \& Sons, Inc., 2012.



\bibitem{Moore} J. Moore, The singular solutions to a singular quadratic minimization problem. International Journal of Control, 1974, 20(3): 383-393.



\bibitem{pinv} R. Penrose, A generalized inverse of matrices, Mathematical Proceedings of the Cambridge Philosophical Society, 52: 17-19, 1955.

\bibitem{mp} L. S. Pontryagin, V. G. Boltyanskii, R. V. Gamkrelidze, E. F. Mishchenko. \emph{The mathematical theory of optimal process}.
English translation. Interscience, 1962.


\bibitem{SCIS2}	Q. Qi, H. Zhang, Time-inconsistent stochastic linear quadratic control for discrete-time systems, \emph{SCIENCE CHINA Information Sciences}, 60(12), 120204:1-120204:13, 2017.

\bibitem{rami} M. Ait Rami, X. Chen, X. Y. Zhou, Discrete-time indefinite LQ control with state and control dependent noise, Journal of Global Optimization, 23: 245-265, 2002.


\bibitem{zhouxunyu} M. A. Rami, J. B. Moore, and X. Y. Zhou, Indefinite stochastic linear
quadratic control and generalized Riccati equation, \emph{SIAM Journal on Control and Optimization}, 40, 4: 1296¨C1311, 2001.



\bibitem{SCIS3}	J. Shi, G. Wang, J. Xiong, Linear-quadratic stochastic Stackelberg differential game with asymmetric information,
\emph{SCIENCE CHINA Information Sciences}, 60, 092202:1-092202:15, 2017.


\bibitem{Speyer} J. Speyer, D. Jacobson, Necessary and Sufficient Conditions for Optimality for Singular Control Problems, \emph{Journal of Mathematical Analysis and Applications}, Vol 31, No. 1, 1971.



\bibitem{sunliyong} J. Sun, X. Li, J. Yong, Open-loop and closed-loop solvabilities for stochastic linear quadratic optimal control problems, \emph{SIAM Journal on Control and Optimization}, 54(5): 2274-2308, 2016.


\bibitem{Williems} J. C. Willems, A. Kitapci, L. M. Silverman, Singular optimal
control: a geometric approach. IAM Journal of Control and Optimization, 1986, 24(2): 323-337.


\bibitem{SCIS1} J. Xu, J. Shi, H. Zhang, A leader-follower stochastic linear quadratic differential game with time delay, \emph{SCIENCE CHINA Information Sciences},
61:112202, 2018.

\bibitem{hszhang} H. Zhang, L. Lin, J. Xu, M. Fu, Linear quadratic regulation and stabilization of
discrete-time Systems with delay and multiplicative noise, \emph{IEEE Transactions on Automatic Control}, 60(10): 2599-2613, 2015.


\bibitem{hszhang1} H. Zhang, J. Xu, Control for It\^{o} stochastic systems with input delay, \emph{IEEE Transactions on Automatic Control}, 62(1): 350-365, 2017.


\bibitem{zhangIR} H. Zhang, J. Xu, Optimal control with irregular performance, \emph{SCIENCE CHINA Information Sciences},
62:192203, 2019.




\bibitem{hszhang1} H. Zhang, J. Xu, Control for It\^{o} stochastic systems with input delay, \emph{IEEE Trans.\ Autom.\ Control}, 62(1): 350-365, 2017.





\bibitem{XuZhang} H. Zhang, X. Zhang, Pointwise second-order necessary conditions for stochastic optimal controls, Part I: The case of convex control constraint.
SIAM Journal of Control and Optimization, 2015, 53(4): 2267-2296.




\end{thebibliography}
\end{document}